\documentclass{paper}
\usepackage[centertags]{amsmath}
\usepackage{amsfonts}
\usepackage{amssymb}
\usepackage[T1]{fontenc}
\usepackage{latexsym}
\usepackage{graphicx}
	
\title{A simplified method for the evaluation of the total resistance of a foiling yacht in upright condition}
\author{Daniele Peri}
\date{Istituto per le Applicazioni del Calcolo "M. Picone", Consiglio Nazionale delle Ricerche, Via dei Taurini, 19 - 00185 Roma - Italy}

\begin{document}

\maketitle

\begin{abstract}
An extremely schematic model of the forces acting an a sailing yacht equipped with a system of foils is here presented and discussed. The role of the foils is to raise the hull from the water in order to reduce the total resistance and then increase the speed. Some CFD simulations are providing the total resistance of the bare hull at some values of speed and displacement, as well as the characteristics (drag and lift coefficients) of the 2D foil sections used for the appendages. A parametric study has been performed for the characterization of a foil of finite dimensions. The equilibrium of the vertical forces and longitudinal moments, as well as a reduced displacement, is obtained by controlling the pitch angle of the foils. The value of the total resistance of the yacht with foils is then compared with the case without foils, evidencing the speed regime where an advantage is obtained, if any.
\end{abstract}

\section{Introduction}

Today, the use of foils in racing yachts has become increasingly common, and this type of yacht seems to be very attractive even for uncompetitive users. The increased speed achieved using foils, along with the intriguing idea of driving a "sailing yacht", makes this solution very popular, but their design is not applicable to every ship and in any condition. Among other things, the presence of these appendices requires an adequate control system, the management of which is not trivial and above all safe for an uneducated user without adequate preparation: Sudden changes in wind/course conditions can cause strong boat reactions, and the likelihood of involvement in dangerous situations is not so unusual. But above all, since the lifting generated by the foils has a price, in terms of induced and viscous resistance, its usefulness is evident for high speeds and very light hulls, but it is not guaranteed for low speeds and for heavy ships.

In this paper we present a simple method for the preliminary evaluation of the total resistance of a ship with foils in upright condition. The required data are very limited and can be obtained with simple (but reliable) tools with a moderate computational effort. The outcome is the new resistance curve of the yacht, able to highlight the range of speeds where there is a real advantage (if any) in the use of foils.

\section{Schematic representation of the acting forces}

In a very schematic way, we can consider our yacht as equipped with three foils. The first couple, that in the following we will recall as {\em main foil} also if they are two (but symmetric), is responsible for a pure lifting force, needed to raise the hull, in part or completely, out of the water: we assume that this foil is mounted in deep and at the same longitudinal position of the center of gravity of the hull. This force is not producing a longitudinal moment. We have also a second foil, hinged at the lowest part of the rudder, which will keep the yacht in a horizontal position. We are here to disregard the sail forces that tend to cause the hull to bow down. Both the foils can be tuned, changing their angle of attack independently. In figure \ref{fig:Eq} we have the graphical representation of this configuration. Under these hypothesis, the procedure for the computation of forces is the following:
\begin{figure}[htb]
	\centering
	\includegraphics[width=\textwidth]{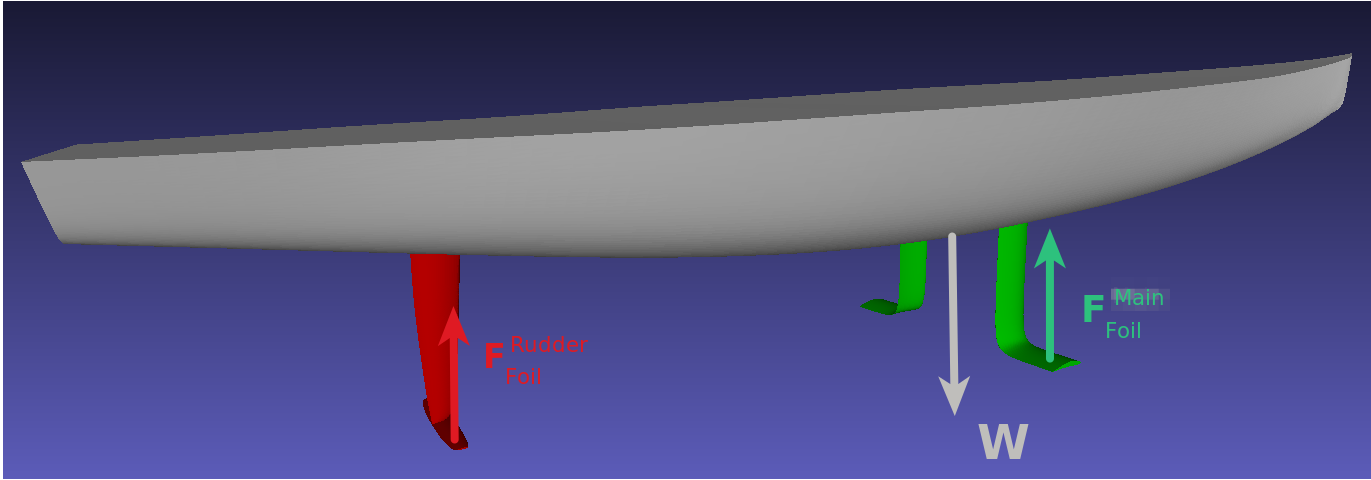}
	\caption{Schematic representation of the forces acting on the hull due to the foils. {\bf W} is the hull weight, the forces generated by the foils are also indicated with {\bf F}.}
	\label{fig:Eq}
\end{figure}

\begin{enumerate}
	\item Fix the total displacement of the yacht.
	\item Fix the displacement of the yacht we want to realize (0 if the yacht is completely out of the water).
	\item Fix the speed of the yacht.
	\item Estimate the resistance of the yacht at the required speed and displacement in horizontal position, without sinkage and trim.
	\item Estimate the forces and moments that tends to change the sinkage and trim of the yacht.
	\item Find the pitch angle of the rudder foil generating the moment able to keep the yacht in equilibrium considering the estimated forces/moment.
	\item Find the pitch angle of the main foil generating the force able to keep the yacht in desired position.
	\item If the estimated pitch angle of the main foil is too high, fix its pitch angle to the maximum allowed value and then find the new displacement of the yacht.
	\item Drag of the foils is computed at the current pitch angles and added to the total resistance.
\end{enumerate}

This procedure requires two main sets of data: one for the hull and one for the foils. There are strong differences between the origins and management of the two datasets: in particular, we can try to apply some strong simplification on the computations of the forces provided by the foils. For this reason, the two sources are here described separately.

\subsection{Hull forces and moments}

For the hull, we need the value of all the forces and moments acting on the hull fixed in horizontal position at a number of speeds and displacements. These data have been obtained by using a CFD solver based on the potential theory \cite{Gadd}. The code has been used excluding the option of the evaluation of the sinkage and trim, so that the vertical force and the longitudinal moment have been computed together with the total resistance. Computations have been performed for three different values of the displacement and for several speeds, in a range of Froude numbers in between 0.20 and 0.65. The viscous effects are not directly included in the adopted mathematical model, so that the results are not fully reliable in the planing and semi-planing regime of the hull. Anyway, for a more accurate use of the overall methodology, viscous computations could be applied instead.

\begin{figure}[htb]
	\centering
	\includegraphics[width=\textwidth]{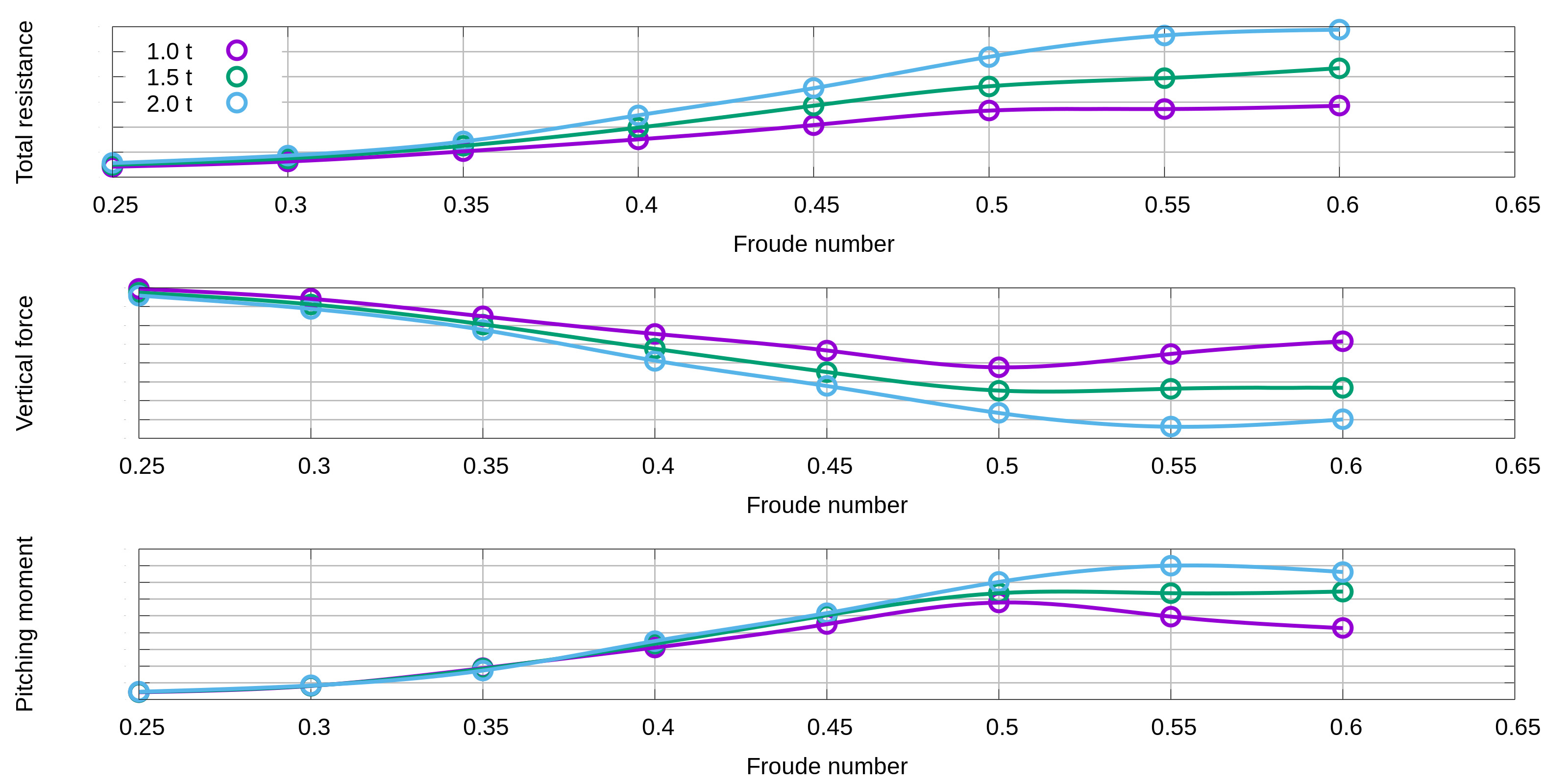}
	\caption{Longitudinal force, vertical force and pitching moment acting on the bare hull in fixed horizontal position in a range of displacements and speeds. Values are hidden due to confidentiality.}
	\label{fig:rt}
\end{figure}

\begin{figure}[htb]
	\centering
	\includegraphics[width=\textwidth]{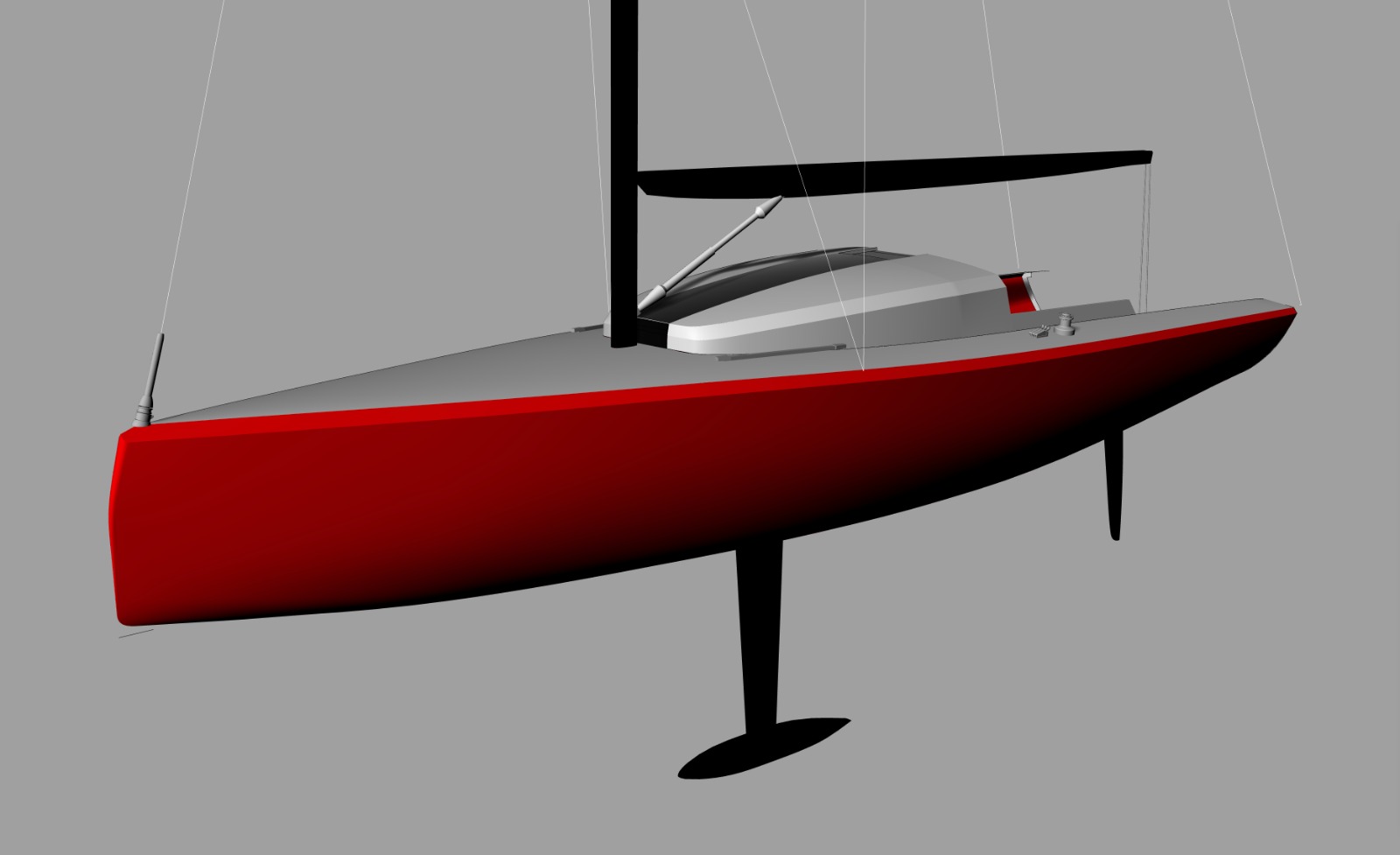}
	\caption{Non-dimensional values of lift and drag as a function of the aspect ratio of the wing (dots). Data are fitted with and exponential-type curve (green line)}
	\label{fig:scafo}
\end{figure}

The adopted hull shape has been provided by the {\em DuckDesign} design team. Due to confidentiality, we cannot show here details of the shape or provide numerical data. A perspective view of the full appended hull is reported in figure \ref{fig:scafo}. The yacht is 9.15 meters long and has a design displacement of 1.5 tons. The overall shape is quite conventional, with a straight bow and a rounded transom stern, not particularly wide.

\begin{figure}[htb]
\centering
\includegraphics[width=\textwidth]{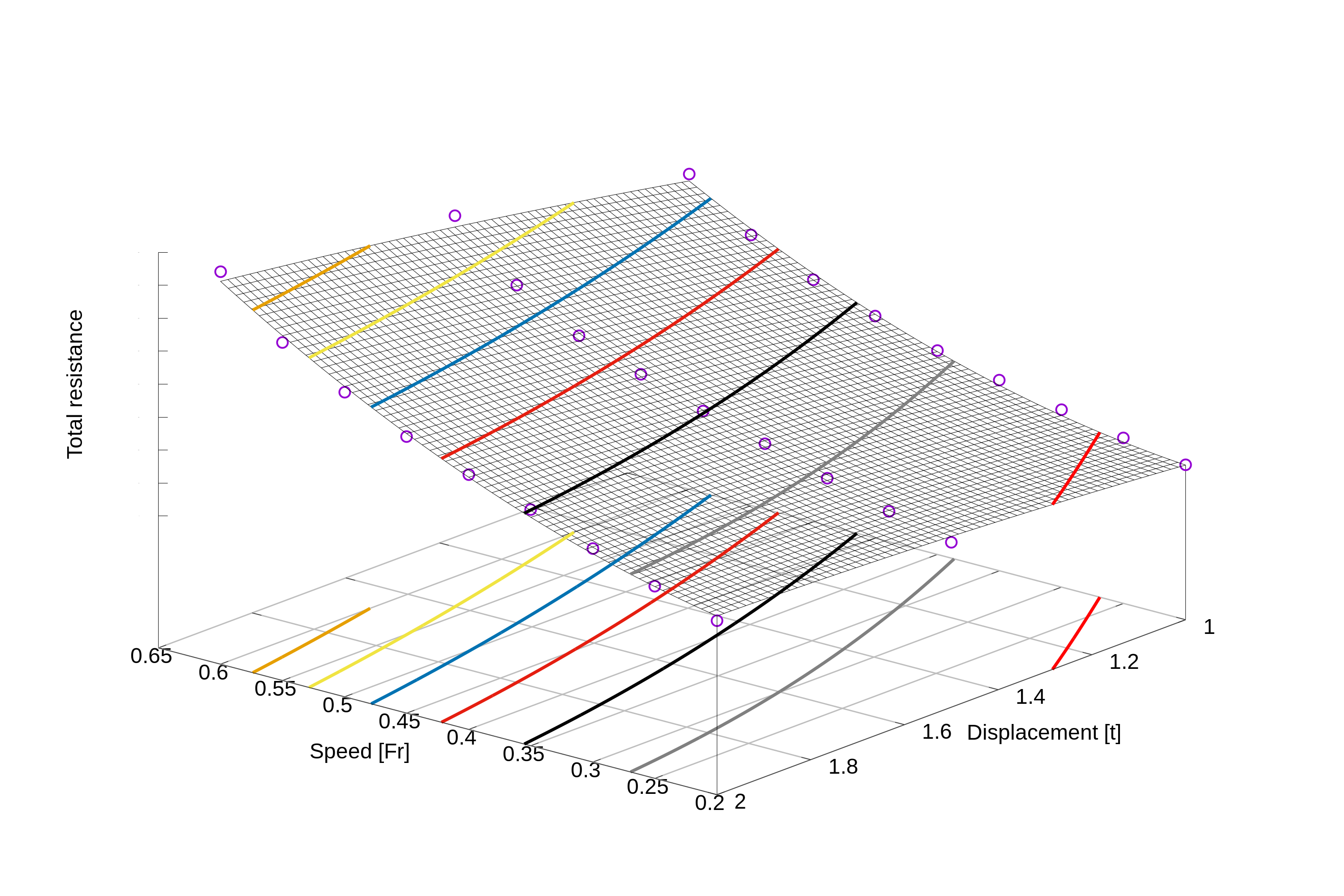}
\caption{Response surface of the total resistance as a function of speed and displacement. Second order polynomial model.}
\label{fig:fit}
\end{figure}

In order to be able to obtain the value of the different forces and moments for every possible value of speed and displacement inside the investigated range, a second order polynomial response surface has been obtained by fitting the available data. Results are reported in figure \ref{fig:fit}.

\subsection{Foil drag and lift}

\begin{figure}[htb]
	\centering
	\includegraphics[width=\textwidth]{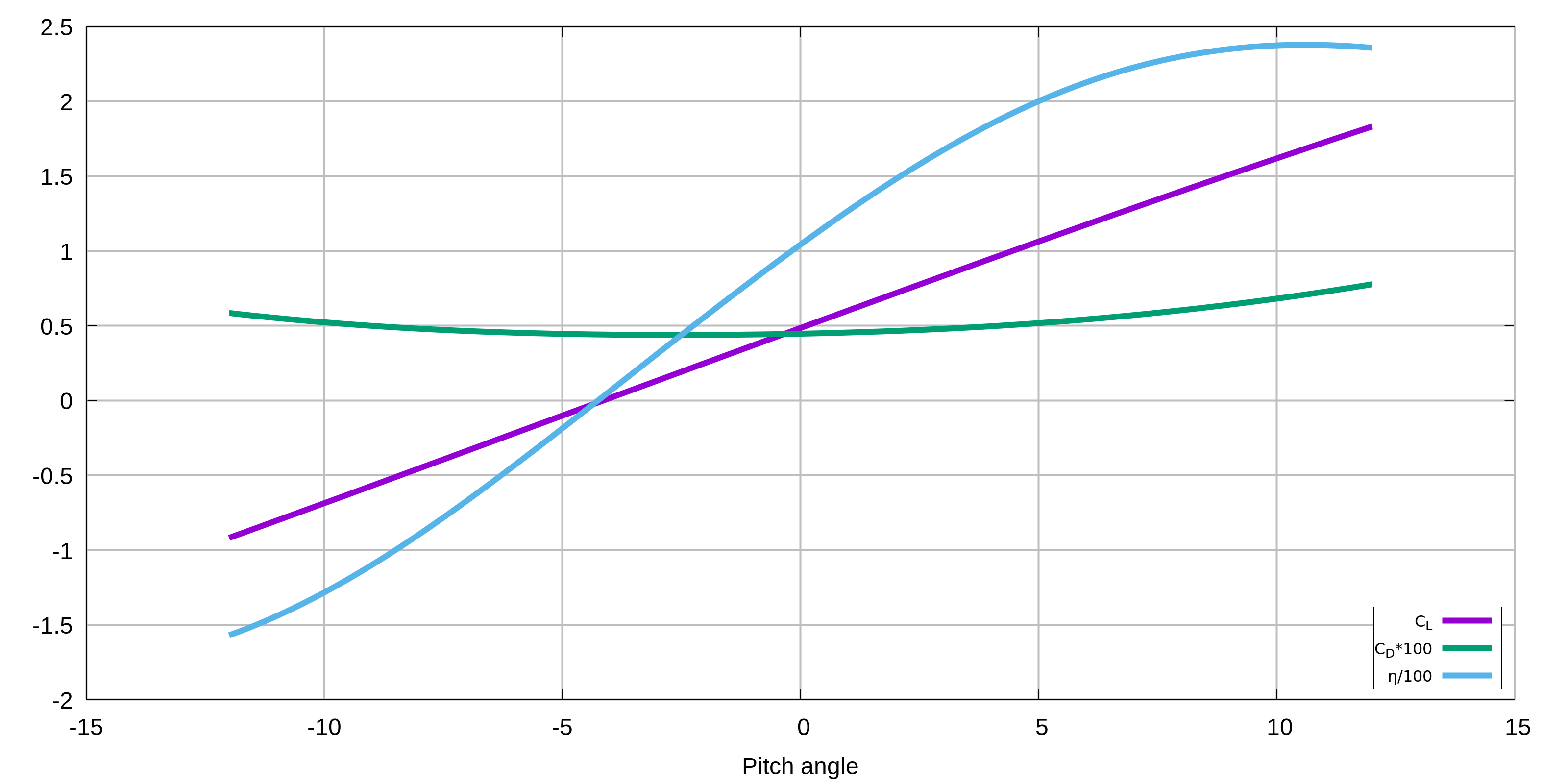}
	\caption{Lift coefficient, drag coefficient and efficiency as a function of the angle of attack for the 2D foil. Data are scaled in order to be of the same order of magnitude on the plot.}
	\label{fig:pol}
\end{figure}

A very easy way to obtain the characteristics of a foil section is the use of a 2D CFD solver. One of the most popular is called {\tt Xfoil} \cite{Xfoil}, a simulation tool including also a simplified (but reliable) formulation for the estimation of the thickness of the boundary layer, so that also the viscous effects, including separation for not extreme geometries, can be evaluated. The solver includes also the generation of the geometry of a wide variety of NACA foils, but a specific foil shape can be also provided as a ordinate sequence of dots. A typical outcome is reported in figure \ref{fig:pol}, where the lift and drag coefficients are plotted together with the efficiency of the foil (defined as the ratio between lift and drag, that is, the amount of lift for unit drag). Data are scaled in order to be of the same order of magnitude on the plot.

Unfortunately, these data represent the ideal performances of the foil, computed in the case of an infinitely wide wing, whereas our foils have finite dimensions. In the case of a finite span, the performances of the foil deteriorates: this is due to the cross-flow generated at the end of the wing, that tends to reduce the pressure difference between the two sides of the wing. There is not a clear relationship with the span. As a consequence, we need to derive this relationship. If we consider the so called {\em aspect ratio} (AR) of the wing, defined as the ratio between span and chord, a 2D computation represents the case of $AR = \infty$. Our relationship could be a function of the AR of the foil.

To this aim, a series of simulations for a number of wings with the same foil section but different values of AR have been performed. The foil section is the NACA 4412, with an angle of attack of 4 degrees. Simulations have been performed by using a CFD solver from the library {\tt OpenFOAM}, in particular the standard unsteady RANSE solver {\tt pimpleFoam}, adopting the realizable-$k\varepsilon$ turbulence model. In the simulations, the wing has a chord of 1 meter and moves at a speed of $20m/s$ (Reynolds number is $2\cdot10^7$,). AR has been varied between 2 and 48.

Using three of the last values of the AR (12, 24 and 48), a convergence study of the lift and drag coefficients as a function of AR has been performed: the asymptotic value thus obtained is theoretically equivalent to the value of a 2D simulation. Results are reported in table \ref{tab:conv}.

\begin{table}
\begin{tabular}{|l|l|l|l|l|l|} \hline
       & AR 12    & AR 24    & AR 48    & Asymptotic & Convergence order \\ \hline
 Drag  &  3034.37 &  2539.55 &  2175.41 &  2054.03   & 0.44              \\ \hline
 Lift  & 80295.58 & 85935.75 & 89683.19 & 90932.33   & 0.59              \\ \hline
\end{tabular}\caption{Convergence study for the dependence of drag and lift from AR of a finite-span wing. For simplicity of representation, here the forces have been divided only by the wing surface ($F_x/S$ and $F_z/S$ are reported)}\label{tab:conv}
\end{table}

\begin{figure}[htb]
	\centering
	\includegraphics[width=\textwidth]{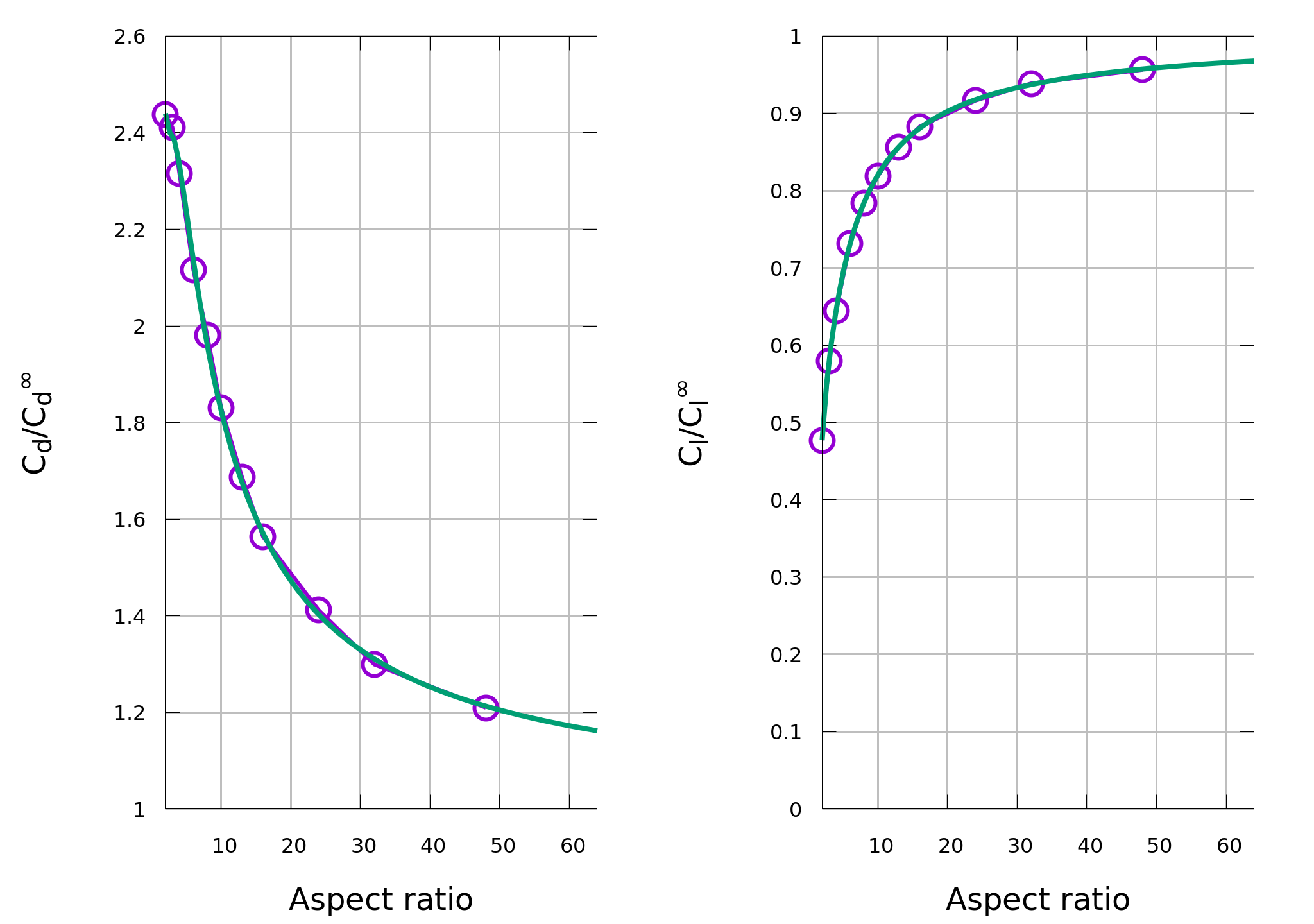}
	\caption{Non-dimensional values of lift and drag as a function of the aspect ratio of the wing (dots). Data are fitted with and exponential-type curve (green line)}
	\label{fig:convAR}
\end{figure}

The computed values of lift and drag coefficients have been divided by the proper asymptotic value, and results are reported in the figure \ref{fig:convAR}. The coefficients represent the values by which it is necessary to multiply the coefficients obtained by the $AR=\infty$ wing when a finite value of AR is to be considered. The data are well aligned with a polynomial curve. By definition, the asymptotic vale of the polynomial curve is 1. In details:

\begin{eqnarray}\label{Cd}
	C_d/C_d^\infty & = &       1.0        -     0.110/AR   \nonumber \\
	               &   & +   217.261/AR^2 -  2742.862/AR^3 \nonumber \\
                   &   & + 16343.289/AR^4 - 51940.417/AR^5 \nonumber \\
                   &   & + 83982.370/AR^6
\end{eqnarray}
\begin{eqnarray}\label{Cl}
	C_l/C_l^\infty & = &       1.0       -      0.239/AR   \nonumber \\
                   &   & -    35.349/AR^2 +   340.353/AR^3 \nonumber \\
                   &   & -  1461.250/AR^4 +  2927.150/AR^5 \nonumber \\
                   &   & -  2191.410/AR^6
\end{eqnarray}

These two relationships can be easily applied in order to correct the 2D data, e.g. obtained by using \cite{Xfoil}, once the main dimensions of the foils are known. In the computations, a reduction of the efficiency of 10\% has been introduced in order to take into account the interference of the devices on which the foils are attached.

A general comment on the fluid dynamic data shown in figure ref{fig:convar}. Unexpectedly, a clear convergence is not yet reached for a fairly large value of the AR: before performing the calculations, it was expected to stop at a much smaller AR value, and also to observe a stabilization of aerodynamic coefficients. On the contrary, especially for the drag coefficient, the convergence is really slow, and for a wing of $AR=48$ we still have an increase in the resistance of about 20\% with respect to the asymptotic value. The lift coefficient is slightly more in line with expectations, and the loss is only about 5\% in the same conditions. From the computational standpoint, the increase of AR also increases (almost linearly) the grid size for the CFD solver (trivially, because the wing is larger), also increasing the computational cost of the simulation: for this reason the study stops at $AR=48$. It should be considered that, since the common values of the AR in yacht design are significantly lower than 48, and the purpose of this study is to derive a reliable relationship for this specific application, it makes no sense to consider even greater AR values, since the calculated values are absolutely sufficient to cover the range of our interest. But a dedicated study would be helpful in the aerodynamic field in general.

\section{Application to a test case}

In the following example, in order to be more adherent to the performed computations, we are using the same foil section for both the foils. The selected foil is obviously the NACA 4412, since all the previous computations for the tuning on the AR have been made on this foil. A deeper investigation would be needed in order to understand the possibility to extend the equations \ref{Cd} and \ref{Cl} for any foil.

The first configuration is the following: a main foil with a chord of 20 cm and a span of 60 cm (dimensions are for the single foil), and rudder foil with a chord of 5 cm and a overall span of 40 cm, 20 for each side of the rudder.

\begin{figure}[htb]
	\centering
	\includegraphics[width=\textwidth]{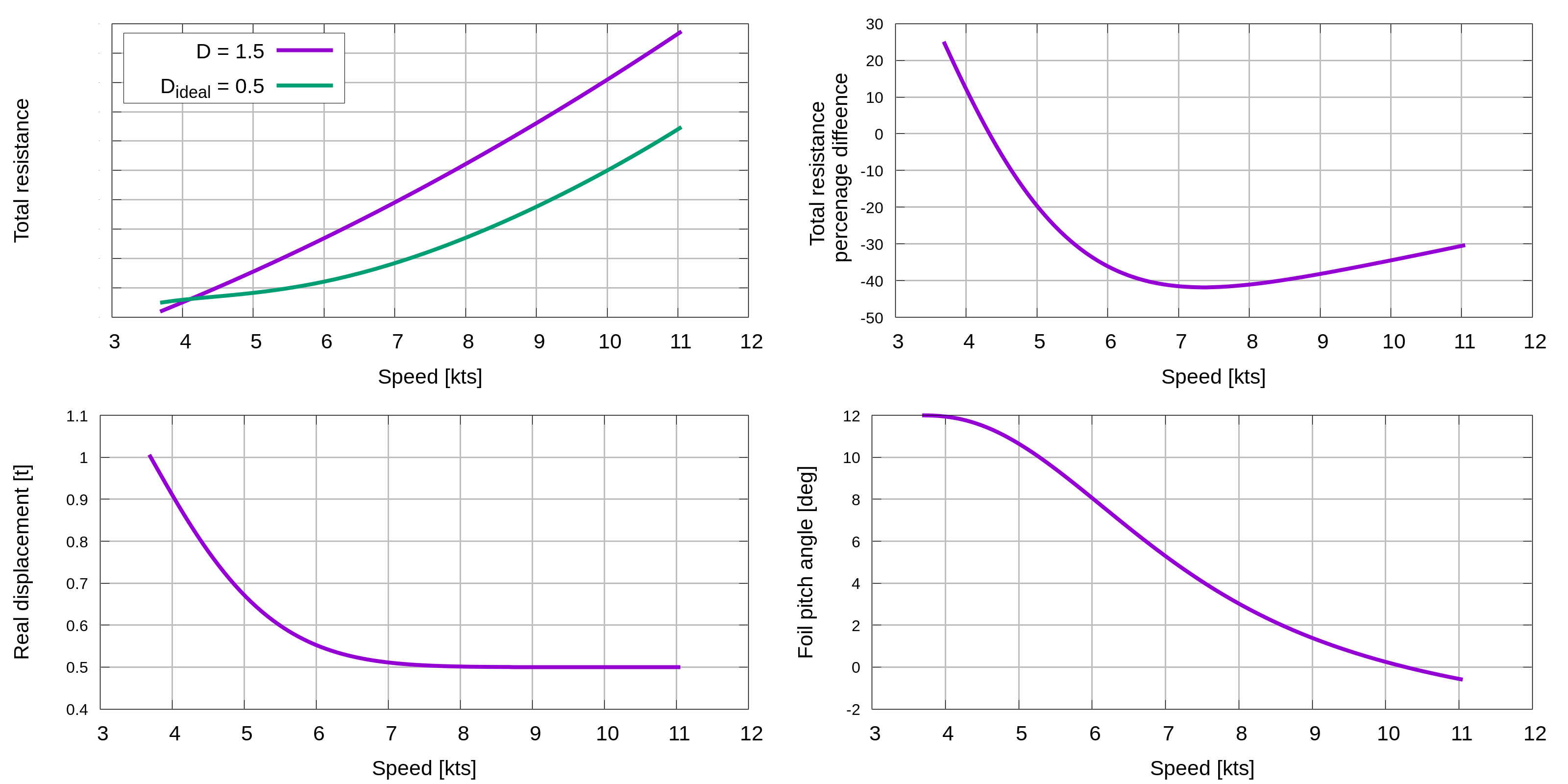}
	\caption{On top, left, the comparison of the total resistance for the hull with and without foils in a range of speeds. Numerical values are hidden due to confidentiality. On top, right, the percentage differences. On bottom, the real displacement of the hull (left) and the angle of attack of the main foil (right) as a function of speed. Main foil has a span of 60 cm and a chord of 20 cm. The rudder foil has a span of 40 cm and a chord of 5 cm.}
	\label{fig:Results}
\end{figure}

Results are reported in figure \ref{fig:Results}. We can see how for a speed lower than 4.5 knots there is no advantage in using the foils, and the presence of the appendages causes the overall performance of the vessel to deteriorate, although the hull is clearly raising up from the water (the weight is reduced to 0.5 tons). On the contrary, for speeds higher than 4.5 knots, the advantage is clear. We have a maximum advantage of about 45\% at the speed of 7 knots. 
  
\begin{figure}[htb]
	\centering
	\includegraphics[width=\textwidth]{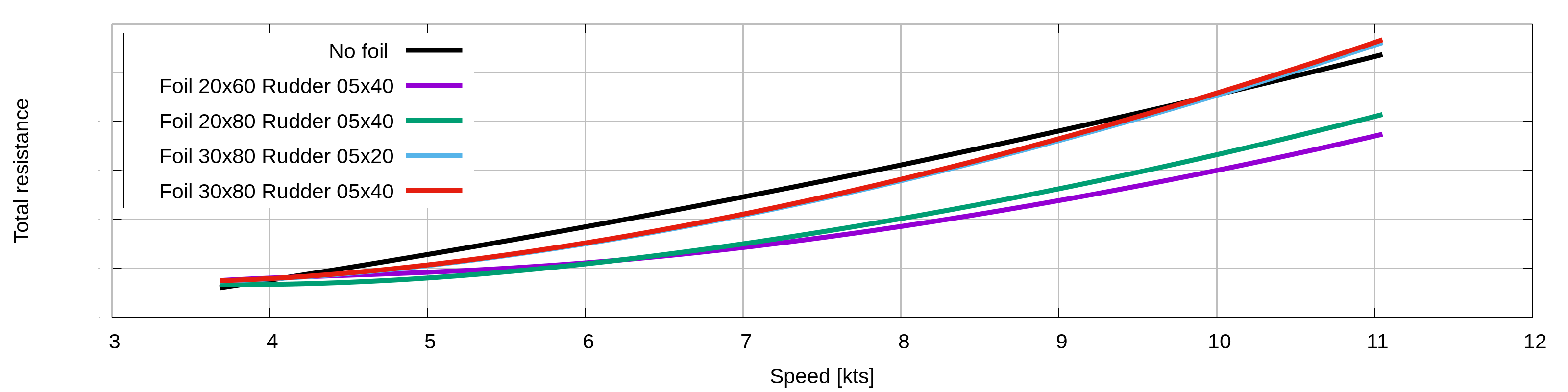}
	\caption{Comparison of the predicted performances for different choices of foils dimensions.}
	\label{fig:Comparison}
\end{figure}

In order to evaluate the different performances with different foil sets, we have checked the effect of a change of the foil dimensions. In figure \ref{fig:Comparison}, 4 different foil configurations are compared in terms of total resistance. Here the increase of the span of the main foil from 60 to 80 cm is further reducing the total resistance in the range of speeds from 4 to 6 knots, but the performances of the previous main foil are better for higher speeds. A possible explanation is that the greater foil surface is able to provide an higher lift force at low speed: the angle of attack is limited in order to prevent separation, and the lifting force is proportional to the foil surface and to the speed squared, so the only way to increase the lift force if the angle of attack is fixed is to increase the surface. But as soon as the lift requirements become less extreme due to the increased speed, and we can apply a smaller angle of attack, the greater surface represents a penalty. On the contrary, an increase of the chord of the main foil is not beneficial. Clearly, this wold be a good application for an optimization algorithm, able to detect the best compromise once the target performances are clearly fixed. Anyway, we can deduce a first general indication from this comparison: it appears that the surface of the appendices must be kept minimal if the boat is very fast, while it must be greater in case of low speeds. 

\section{Conclusions}

A simple tool for the preliminary evaluation of the performances of a sailing yacht equipped by a system of foils is presented. Although numerical results are not validated, the tool proves to be able to give indications about the different performances provided by different set ups.

The verification of the results obtained, e.g. using high quality CFD simulations of the full yacht, could allow the verification of the results.

Another point is represented by the study on the effects of the aspect ratio on the performances of a wing of finite dimensions: the very preliminary study here presented should be extended to a greater variety of foils.

\section*{Acknowledgments}

The author is grateful to Arch. Massimo Paperini for the hull design.


\end{document}